\documentclass[11pt]{amsart}

\usepackage[english]{babel} 
\usepackage{graphicx}
\usepackage{amsmath,amsfonts,amssymb}
\usepackage{cite}
\usepackage{color}  
\usepackage{float}
\usepackage{caption}
\usepackage{multirow}
\usepackage{verbatim}
\usepackage{hyperref}
\usepackage{xcolor}
\hypersetup{
	colorlinks,
	linkcolor={red!50!black},
	citecolor={blue!50!black},
	urlcolor={blue!80!black}
}
\usepackage{enumerate}
\newcommand{\sprod}[2]{\left\langle #1, #2\right\rangle}
\newcommand{\R}{{\mathbb R}}
\newcommand{\N}{{\mathbb N}}
\renewcommand{\v}[1]{\ensuremath{\mathbf{#1}}}
\renewcommand{\b}{\v{b}}
\newcommand{\x}{\v{x}}
\newcommand{\y}{\v{y}}
\newcommand{\w}{\v{w}}
\newcommand{\z}{\v{z}}
\newcommand{\lb}{\boldsymbol{\lambda}}

\newcommand{\bd}{\v{b}^\delta}
\newcommand{\ts}[1]{\mbox{\textsf{#1}}}
\newcommand{\norm}[1]{\left\| #1 \right\|}

\date{}
\usepackage{ifthen}
\usepackage{float}
\floatstyle{ruled}
\newfloat{algorithm}{htb}{alg}
\floatname{algorithm}{Algorithm}
\newcounter{algo@row}
\newcounter{algo@rowindent}
\newcommand{\algofont}[1]{\textbf{#1}}
\newcommand{\algonumbersize}[1]{\scriptsize{#1}}
\newcommand{\algopreitem}[1][\arabic{algo@row}]{\texttt{\algonumbersize{#1}}}
\newcommand{\algoitemskip}{\hspace{\value{algo@rowindent}cc}}
\newenvironment{algo}{\vskip.3em\small%
	\begin{list}{\algopreitem\texttt{\algonumbersize{:}}}{%
			\usecounter{algo@row}%
			\setcounter{algo@rowindent}{0}%
			\setlength{\itemindent}{2em}%
			\setlength{\labelwidth}{2em}
			\setlength{\parsep}{0cm}%
		}%
	}{
	\end{list}
}

\newtheorem{theorem}{Theorem}[section]

\newtheorem{remark}[theorem]{Remark}

\numberwithin{equation}{section}
\numberwithin{theorem}{section}
\newcommand{\algonewnestedopen}[2]{
	\newcommand{#1}[1][]{%
		\ifthenelse{\equal{##1}{}}{\item}{\item[{\algopreitem[##1]}]}
		\algoitemskip\algofont{#2}%
		\addtocounter{algo@rowindent}{1}%
		\ignorespaces
	}
}
\newcommand{\algonewnestedaux}[2]{
	\newcommand{#1}[1][]{
		\addtocounter{algo@rowindent}{-1}
		\ifthenelse{\equal{##1}{}}{\item}{\item[{\algopreitem[##1]}]}
		\algoitemskip\algofont{#2}%
		\addtocounter{algo@rowindent}{+1}%
		\ignorespaces
	}
}
\newcommand{\algonewnestedclose}[2]{
	\newcommand{#1}[1][]{
		\addtocounter{algo@rowindent}{-1}
		\ifthenelse{\equal{##1}{}}{\item}{\item[{\algopreitem[##1]}]}
		\algoitemskip\algofont{#2}%
		\ignorespaces
	}
}
\newcommand{\algonewcommand}[2]{
	\newcommand{#1}[1][default]{
		\ifthenelse{\equal{##1}{default}}{\item}{\item[{\algopreitem[##1]}]}%
		\algoitemskip\algofont{#2}%
		\ignorespaces
	}%
}
\newcommand{\algonewkeyword}[2]{\newcommand{#1}{\algofont{#2}}}
\algonewcommand{\STATE}{\ignorespaces}
\algonewcommand{\INPUT}{Input: }
\algonewcommand{\pINPUT}{\phantom{Input: }}
\algonewcommand{\COMPUTE}{Compute: }
\algonewcommand{\OUTPUT}{Output: }
\algonewcommand{\pOUTPUT}{\phantom{Output: }}
\algonewnestedopen{\IF}{if }
\algonewnestedaux{\ELSEIF}{else if }
\algonewnestedaux{\ELSE}{else }
\algonewnestedclose{\ENDIF}{end if }
\algonewnestedopen{\FOR}{for }
\algonewnestedclose{\ENDFOR}{end for }
\algonewnestedopen{\PARFOR}{parfor }
\algonewnestedclose{\ENDPARFOR}{end parfor }
\algonewnestedopen{\WHILE}{while }
\algonewnestedclose{\ENDWHILE}{end while }
\algonewcommand{\BREAK}{break}%
\algonewkeyword{\For}{for }%
\algonewkeyword{\To}{to }%
\algonewkeyword{\Do}{do }%
\algonewkeyword{\If}{if }%
\algonewkeyword{\Then}{then }%
\algonewkeyword{\Else}{else }%
\algonewkeyword{\End}{end }%
\algonewkeyword{\AND}{and }%
\algonewkeyword{\True}{true }%
\algonewkeyword{\False}{false }%
\algonewkeyword{\Call}{call }%
\algonewkeyword{\irbleigs}{irbleigs }%
\algonewkeyword{\tridiag}{tridiag}%
\algonewkeyword{\reorth}{reorth}%

\usepackage[normalem]{ulem}
\newcommand{\stkout}[1]{\ifmmode\text{\sout{\ensuremath{#1}}}\else\sout{#1}\fi}
\usepackage{todonotes}

\begin{document}
\title{Graph Laplacian for image deblurring\footnote{This is a preprint.}}
\author{Davide Bianchi}
\address{Dipartimento di Scienze e Alta Tecnologia\\
	Universit\`a dell'Insubria\\
	via Valleggio 11\\
	I-22100 Como, Italy}
\email{d.bianchi9@uninsubria.it}

\author{Alessandro Buccini}
\address{Department of Mathematics and Computer Science\\
	University of Cagliari\\
	09124 Cagliari, Italy}
\email{alessandro.buccini@unica.it}

\author{Marco Donatelli}
\address{Dipartimento di Scienze e Alta Tecnologia\\
	Universit\`a dell'Insubria\\
	via Valleggio 11\\
	I-22100 Como, Italy}
\email{marco.donatelli@uninsubria.it}
\author{Emma Randazzo}
\address{Dipartimento di Scienze e Alta Tecnologia\\
	Universit\`a dell'Insubria\\
	via Valleggio 11\\
	I-22100 Como, Italy}
\email{e.randazzo@uninsubria.it}
\maketitle

\begin{abstract}
Image deblurring is relevant in many fields of science and engineering. To solve this problem, many different approaches have been proposed and among the various methods, variational ones are extremely popular. These approaches are characterized by substituting the original problem with a minimization one where the functional is composed of two terms, a data fidelity term and a regularization term. In this paper we propose, in the classical $\ell^2-\ell^1$ minimization with the non-negativity constraint of the solution, the use of the graph Laplacian as regularization operator. Firstly, we describe how to construct the graph Laplacian from the observed noisy and blurred image. Once the graph Laplacian has been built, we 
solve efficiently the proposed minimization problem splitting the convolution operator and the graph Laplacian by the alternating direction method of multipliers (ADMM). Some selected numerical examples show the good performances of the proposed algorithm.
\end{abstract}
	

\section{Introduction}
We are concerned with the solution of the image deblurring problem; see, e.g., \cite{HNO} for more details on image deblurring. We will assume that the blurring is space-invariant obtaining a linear system of equations
$$
A\x=\b,
$$
where $\x\in\R^N$ and $\b\in\R^M$ are samplings of the unknown image to recover and the blurred image, respectively, while $A\in\R^{M\times N}$ is a structured matrix (see below) whose singular values decay rapidly and with no significant gap. The discretization process, along with measurement errors, introduces some perturbations in the data, namely $\boldsymbol{\eta}\in\R^M$ such that $\| \b-\boldsymbol{\eta} \|=\delta$ in the Euclidean norm, leading to the system
$$
A\x = \b + \boldsymbol{\eta}=\b^\delta.
$$
The perturbations $\boldsymbol{\eta}$ are often referred to as noise. Since, in general, $\boldsymbol{\eta}\notin\mathcal{R}(A)$, we would like to solve the least-square problem
\begin{equation}\label{eq:ls}
\arg\min_\x\norm{A\x-\b^\delta},
\end{equation}
where $\norm{\cdot}$ is the Euclidean norm. Let $A^\dagger$ denote the Moore-Penrose pseudo-inverse of $A$. The naive solution of \eqref{eq:ls}, $A^\dagger\b^\delta$, is usually a poor approximation of the desired solution $\x^\dagger=A^\dagger\b$; see, e.g., \cite{H96,H98} for more details. This is due to the fact that $A$ is severely ill-conditioned and the solution of \eqref{eq:ls} is extremely sensitive to the presence of the noise in the data. To compute an approximation of $\x^\dagger$ we need to use regularization methods. These methods aim at reducing the sensitivity mentioned before. One of the most popular regularization method is Tikhonov regularization where the original minimization problem is substituted by another one of the form
\begin{equation}\label{eq:Tikhnov}
\arg\min_\x\norm{A\x-\b^\delta}^2+\mu\norm{L\x}^2,
\end{equation}
where $L\in\R^{p\times n}$ is an operator such that $\mathcal{N}(A)\cap\mathcal{N}(L)=\{\v{0}\}$. The matrix $L$ is the so-called regularization operator and its role is to enforce some a-priori knowledge on the reconstruction. If $A$ is the discretization of an integral operator, then usually $L$ is chosen to be a discretization of either the first or the second derivative; see, e.g., \cite{DR14}. The formulation \eqref{eq:Tikhnov} can be extended for a general $\ell_p$-norm
\begin{equation}\label{eq:l2lp}
\arg\min_\x\norm{A\x-\b^\delta}^2+\mu\norm{L\x}^p_p,
\end{equation}
where $\norm{\x}_p^p=\sum_{i=1}^{n}|x_i|^p$ for $p>0$. Note that, for $p<1$ the function $\x\mapsto\norm{\x}_p$ is not a norm; see, e.g., \cite{HLMSR17,LMSR15,BR19,BRXXb,BPR20, EGLT17,CH14}. In this paper, we consider le graph Laplacian for $L$ in \eqref{eq:l2lp} with $p=1$. Therefore, our minimization problem is of the form
\begin{equation}\label{eq:constr_l2l1}
\arg\min_{\x\geq0}\norm{A\x-\b^\delta}^2+\mu\norm{L\x}_1,
\end{equation}
where we have introduced the non-negativity constraint since images cannot attain negative values.

Recently, the graph Laplacian of a network, built from the image itself, has been proposed as regularization operator mainly for image denoising; see. e.g., \cite{MS14,KM15,PC17,KF09,LN19,YAM16,SPKV15}. In this paper, we build such an operator and use it in \eqref{eq:constr_l2l1} to reconstruct blurred images. We propose an automated procedure for the construction of an appropriate graph Laplacian and show its performances in some selected numerical examples. We compare the graph Laplacian with the standard Total Variation (TV) operator (see \cite{ROF92}) and show that our proposal can lead to substantial improvements in the quality of the reconstructed images.

This paper is structured as follows: in Section~\ref{sect:Constr} we recall the definition of the Laplacian of a given graph and we construct the one that we use in the following. Section~\ref{sect:deblur} presents our algorithmic proposal for the solution of \eqref{eq:constr_l2l1} and Section~\ref{sect:num} contains some selected numerical experiments. Finally, we draw some conclusions in Section~\ref{sect:concl}.

\subsection{Notation}
Discretized images are made by union of several pixels in the plane, and therefore are well represented by nonnegative two-dimensional discrete functions $\x:\R^{n_1\times n_2}\to [0,+\infty)$, $\x(i_1,i_2) = x_{i_1,i_2}\in \R_+$ for $i_1=1,\ldots, n_1,\, i_2=1,\ldots, n_2$. Clearly, this choice is not unique. Since all the operations and analysis we will produce are invariant with respect to the ordering, for the sake of notational simplicity, we fix $n_1=n_2=n=\sqrt{N}$ and consider the lexicographic one-dimensional ordering, that is, $(i_1,i_2)=:i<j:=(j_1,j_2)$ if $i_1<j_1$ or $i_1=j_1$ and $i_2<j_2$. With this choice, an image reads as $\x : \R^N \to [0,+\infty)$, $\x(i)=x_i \in \R_+$ for $i=1,\ldots,n$, and it is said that $\x$ is the vectorization of a square image.

\section{Construction of the graph Laplacian}\label{sect:Constr}
In this section, we first describe how to construct the Laplacian of a given weighted graph. Then we show how to build an appropriate graph, i.e., a graph whose Laplacian is a ``good'' regularization operator, given a good approximation of the exact image $\x^\dagger$. Finally, we provide an algorithm to construct our $L$ given the problem \eqref{eq:ls}.

Given a countable measure space $\left(\mathcal{V},\nu\right)$, where $\nu$ is a positive measure, a symmetric non-negative function $\omega : \mathcal{V}\times \mathcal{V}\to [0,+\infty)$ with zero diagonal is called undirected \emph{graph} on $\mathcal{V}$. The elements $i,j$ of the set $\mathcal{V}$ are called \emph{vertices}, and two vertices are \emph{connected} if $\omega(i,j)>0$. The positive value $w(i,j)$ is called \emph{weight} associated to the \emph{edge} $\{i,j\}$; for a modern analytic introduction to graph theory we refer to \cite{Radek2021}. If $\mathcal{V}$ is a finite set of $n$ elements, then the graph $\omega$ can be uniquely represented, unless permutations, by the \emph{adjacency matrix} $\Omega \in \R^{N\times N}$
$$
\left(\Omega\right)_{i,j}:= \omega(i,j).
$$ 
The linear operator $L_\omega: C(\mathcal{V}) \to C(\mathcal{V})$, acting on the space $C(\mathcal{V}):= \left\{\x : \mathcal{V} \to \R \right\} \simeq \R^N$ via
$$
L_\omega\x(i) := \frac{1}{\nu(i)}\sum_{j} \omega(i,j)\left(\x(i) -\x(j)\right)
$$
is the graph Laplacian on $C(\mathcal{V})$ associated to the graph $\omega$. It is a symmetric operator with respect to the inner product
$$
\langle \x, \y \rangle := \sum_{i} \x(i)\y(i)\nu(i).
$$
In many applications, a quite standard choice for the measure $\nu$ is the \emph{degree} function $\nu=\deg$, defined by
$$
\deg(i):= \sum_{j}\omega(i,j).
$$ 
It measures the whole intensity of the weights associated to the vertex $i$. Clearly, this choice makes $L_\omega$ not symmetric with respect to the standard Euclidean inner product. A good compromise is to choose the homogeneous measure associated to the Frobenius norm of $\Omega$, i.e., $\nu(i)\equiv \|\Omega\|_F$. Let us observe that, writing $D$ as the diagonal matrix such that $(D)_{i,i}=\deg(i)$, then it is easy to check that
$$
\|D\|_1 \leq \|\Omega\|_F \leq \|D\|_2.
$$  
Henceforth, we will assume $\nu=\|\Omega\|_F$. In matrix form, then the graph Laplacian reads
$$
L_\omega = \frac{D-\Omega}{\|\Omega\|_F}.
$$
We wish to construct a graph $\omega$ so that $L_\omega$ can be used in \eqref{eq:constr_l2l1}. In principle we would like to construct $\omega$ such that
$$
\norm{L_\omega\x^\dagger}\approx 0.
$$
To this aim let $\x^*$ be a good approximation of $\x^\dagger$. Define $\omega$ as the weighted and undirected graph on $\mathcal{V}$, whose nodes are the pixels of $\x^\dagger$ and the weights are defined by
$$
\omega(i,j)=\left\{\begin{array}{ll}
{\rm e}^{-(x^*(i)-x^*(j))^2/\sigma}&\mbox{if }i\neq j \mbox{ and }\norm{i-j}_\infty\leq R,\\
0&\mbox{otherwise,}
\end{array}\right.
$$
where $\sigma>0$ and $R\in\N$ are user-defined parameters. Let us recall that we are using the one-dimensional lexicographic ordering, but the nodes of $\mathcal{V}$ represent points in $\R^2$. Therefore, $i=(i_1,i_2), j=(j_1,j_2)$ and $\|i-j\|_\infty=\max\left\{ |i_1-j_1|; |i_2-j_2| \right\}$. Intuitively, the graph is constructed as follows: we connect two pixels if they are close enough and we weight their connection depending on how similar their values are. In particular, we give a strong connection to pixels that have similar values. The parameter $R$ determines how large is the neighborhood we consider for each pixel and $\sigma$ determines how strong the connection between two close pixels should be.

The construction of this graph, and consequently of the graph Laplacian, in turn depends on the construction of an appropriate approximation $\x^*$. As we show in Section \ref{sect:num}, if we could choose $\x^*=\x^\dagger$ we would obtain an almost optimal result. However, this is not possible in realistic scenarios. Therefore, we wish to provide a practical way to determine a good enough $\x^*$ in a totally automatic way. To compute $\x^*$ we propose to solve \eqref{eq:Tikhnov} with a regularization operator defined as follows. Let $L_1\in\R^{n\times n}$ be
$$
L_1=\begin{bmatrix}
-1&1\\
&-1&1\\
&&\ddots&\ddots\\
&&&-1&1\\
1&&&&-1
\end{bmatrix},
$$
i.e., $L_1$ is a discretization of the first derivative with periodic boundary conditions (BCs). Let $I_{n}$ be the identity matrix of order $n$, we define $L_{\rm TV}$ by
\begin{equation}\label{eq:TV}
L_{\rm TV}=\begin{bmatrix}
L_1\otimes I_{n}\\ I_{n}\otimes L_1
\end{bmatrix}\in\R^{2N\times N}.
\end{equation}
Note that $L_{\rm TV}$ is an extremely sparse matrix formed by two Block Circulant with Circulant Blocks (BCCB) matrices stacked one over the other. Therefore, matrix-vector products involving $L_{\rm TV}$ can be performed extremely cheaply (in particular, the flop count is $O(N)$) and $L^T_{\rm TV}L_{\rm TV}$ is a BCCB matrix. We exploit the latter property below.

To simplify the computations we impose periodic BCs to the matrix $A$. Thanks to this choice, $A$ is a BCCB matrix; see \cite{HNO} for more details. We recall that BCCB matrices are diagonalized by the two-dimensional Fourier matrix. Let $F_1\in\R^{n\times n}$ be the Fourier matrix, i.e., $(F_1)_{j,k}={\rm e}^{2\pi\iota (j-1)(k-1)/n}$ with $\iota^2=-1$, then the two-dimensional Fourier matrix is defined by $F=F_1\otimes F_1$. Note that matrix-vector products with $F$ and its inverse $F^*$ can be performed in $O(N\log N)$ flops with the aid of the \texttt{fft} and \texttt{ifft} algorithms.

As discussed above we wish to solve \eqref{eq:Tikhnov} with $L$ described above to determine $L_\omega$, i.e., we wish to solve a problem of the form
\begin{equation}\label{eq:TV_Tik_mu}
	\x_\mu=\arg\min_\x\norm{A\x-\bd}^2+\mu\norm{L_{\rm TV}\x}^2,
\end{equation} 
for a certain $\mu>0$. Thanks to the structure of $A$ and $L_{\rm TV}$ this can be solved cheaply for any $\mu$. We can write
\begin{equation}\label{eq:dec}
A=F^*\Sigma F\quad\mbox{and}\quad L_{\rm TV}=\begin{bmatrix}
	F^*\Lambda_x F\\F^*\Lambda_y F 
\end{bmatrix},
\end{equation}
where $\Sigma$, $\Lambda_x$, and $\Lambda_y$ are diagonal matrices whose diagonal entries are the eigenvalues of $A$, $L_1\otimes I_n$, and $I\otimes L_1$, respectively. We recall that the eigenvalue of a BCCB matrix $C$ can be computed by $F\v{c}_1$, where $\v{c}_1$ is the first column of $C$. Assuming that $\mathcal{N}(A)\cap\mathcal{N}(L_{\rm TV})=\{\v{0}\}$ we have that
\begin{align*}
		\x_\mu&=(A^TA+\mu L^T_{\rm TV}L_{\rm TV})^{-1}A^T\bd\\
		&=\left(F^*\Sigma^* FF^*\Sigma F+\mu \begin{bmatrix}
			F^*\Lambda_x^* F&F^*\Lambda_y^* F 
		\end{bmatrix}\begin{bmatrix}
		F^*\Lambda_x F\\F^*\Lambda_y F 
	\end{bmatrix}\right)^{-1}F^*\Sigma^* F\bd\\
&=\left(F^*\Sigma^*\Sigma F+\mu F^*(\Lambda_x^*\Lambda_x+\Lambda_y^*\Lambda_y)F\right)^{-1}F^*\Sigma^* F\bd\\
&=F^*\left(\Sigma^*\Sigma+\mu(\Lambda_x^*\Lambda_x+\Lambda_y^*\Lambda_y)\right)^{-1}\Sigma^* F\bd,
\end{align*} 
where the matrix to be inverted is a diagonal matrix. Therefore, $\x_\mu$ can be computed for any $\mu$ cheaply.

We now wish to determine in an automatic way the parameter $\mu$. We employ the Generalized Cross Validation (GCV).

Denote by $G(\mu)$ the following function
$$
G(\mu)=\frac{\norm{A\x_\mu-\bd}^2}{\ts{trace}(I-A(A^TA+\mu L^T_{\rm TV}L_{\rm TV})^{-1}A^T)^2}.
$$
The GCV parameter is the minimizer of $G(\mu)$, i.e., $\mu_{\rm GCV}=\arg\min_\mu G(\mu)$.  Given the decomposition \eqref{eq:dec} the value of $G(\mu)$ can be computed in a straightforward way. Introduce the following notation
$$
r_\mu=\norm{A\x_\mu-\bd}\quad\mbox{and}\quad t_\mu=\ts{trace}(I-A(A^TA+\mu L^T_{\rm TV}L_{\rm TV})^{-1}A^T),
$$
therefore, $G(\mu)=r_\mu^2/t_\mu^2$. Using the spectral decomposition of $A$ we have
\begin{align*}
	r_\mu&=\norm{A\x_\mu-\bd}=\norm{F^*\Sigma F F^*\left(\Sigma^*\Sigma+\mu(\Lambda_x^*\Lambda_x+\Lambda_y^*\Lambda_y)\right)^{-1}\Sigma^* F\bd -\bd}\\
	&=\norm{\Sigma\left(\Sigma^*\Sigma+\mu(\Lambda_x^*\Lambda_x+\Lambda_y^*\Lambda_y)\right)^{-1}\Sigma^* F\bd -F\bd}\\
	&=\norm{(\Sigma\Sigma^*\left(\Sigma^*\Sigma+\mu(\Lambda_x^*\Lambda_x+\Lambda_y^*\Lambda_y)\right)^{-1}-I)\widehat{\bd}},
\end{align*}
where $\widehat{\bd}=F\bd$. We now move to the computation of $t_\mu$
 \begin{align*}
 	t_\mu&=\ts{trace}(I-A(A^TA+\mu L^T_{\rm TV}L_{\rm TV})^{-1}A^T)\\
 	&=\ts{trace}(I-F^*\Sigma(\Sigma^*\Sigma+\mu \left(\Lambda_x^*\Lambda_x+\Lambda_y^*\Lambda_y)\right)^{-1}\Sigma^*F)\\
 	&=\ts{trace}(I-\Sigma\Sigma^*(\Sigma^*\Sigma+\mu \left(\Lambda_x^*\Lambda_x+\Lambda_y^*\Lambda_y)\right)^{-1}).
 \end{align*}
We can observe that, once the decompositions \eqref{eq:dec} and $\widehat{\bd}$ have been computed, the evaluation of $G(\mu)$ can be done in $O(N)$ flops. This allows for an extremely fast determination of $\mu_{\rm GCV}$. Finally, we select as $\x^*$ the solution of the minimization problem
\begin{equation}\label{eq:TV_Tik}
\x^*=\arg\min_\x\norm{A\x-\bd}^2+\mu_{\rm GCV}\norm{L_{\rm TV}\x}^2.
\end{equation} 

\begin{remark}
	If $A$ is constructed with BCs different from the periodic ones it is still possible to compute a fairly accurate approximation of $G(\mu)$ using Krylov subspace methods; see \cite{FRR,FRRS,BR21,BXX}.
\end{remark}

We summarize the procedure to construct $L_\omega$ in Algorithm \ref{algo:1}.

\begin{algorithm}\caption{Construction of $L_\omega$ for image deblurring}
	\begin{algo}\label{algo:1}
	\INPUT $A\in\R^{N\times N}$, $\bd\in\R^N$, $R>0$, $\sigma>0$
	\OUTPUT $L_\omega\in\R^{n\times n}$
	\STATE Construct $L_{\rm TV}$ as defined in \eqref{eq:TV}
	\STATE $\Sigma=\ts{diag}\left(F(A_{(:,1)})\right)$
	\STATE $\Lambda_x=\ts{ diag}\left(F((L_{\rm TV})_{(:,1:N)})\right)$
	\STATE $\Lambda_y=\ts{ diag}\left(F((L_{\rm TV})_{(:,N+1:2N)})\right)$
	\STATE $\widehat{\bd}=F\bd$
	\STATE $\mu_{\rm GCV}=\arg\min_\mu \frac{r_\mu^2}{t_\mu^2}$, where $$
	\left\{\begin{array}{l}
		r_\mu=\norm{(\Sigma\Sigma^*\left(\Sigma^*\Sigma+\mu(\Lambda_x^*\Lambda_x+\Lambda_y^*\Lambda_y)\right)^{-1}-I)\widehat{\bd}}\\
		t_\mu=\ts{trace}(I-\Sigma\Sigma^*(\Sigma^*\Sigma+\mu \left(\Lambda_x^*\Lambda_x+\Lambda_y^*\Lambda_y)\right)^{-1})
	\end{array}\right.
$$
	\STATE $\x^*=F^*\left(\Sigma^*\Sigma+\mu(\Lambda_x^*\Lambda_x+\Lambda_y^*\Lambda_y)\right)^{-1}\Sigma^* \widehat{\bd}$
	\STATE Construct $\Omega\in\R^{n\times n}$ as
	$$
	\omega(i,j)=\left\{\begin{array}{ll}
		{\rm e}^{-(x^*(i)-x^*(j))^2/\sigma}&\mbox{if }i\neq j \mbox{ and }\norm{i-j}_\infty\leq R,\\
		0&\mbox{else,}
	\end{array}\right.
	$$
	with $i$ and $j$ representing two-dimensional indexes, sorted in lexicographic order
	\STATE $D=\ts{diag}\{\sum_{j=1}^{n}(\Omega)_{(i,j)}\}$
	\STATE $L_\omega=\frac{D-\Omega}{\norm{\Omega}_F}$
	\end{algo}
\end{algorithm}

\section{Graph Laplacian deblurring}\label{sect:deblur}
We now describe the non-linear model we employ to compute an approximate solution of \eqref{eq:ls}. We consider the graph Laplacian $L_\omega$ constructed by Algorithm \ref{algo:1} and we use it in \eqref{eq:constr_l2l1}.
Therefore, we wish to solve the following
\begin{equation}\label{eq:model}
\arg\min_{\x\geq0}\frac{1}{2}\norm{A\x-\bd}^2+\mu\norm{L_\omega\x}_1.
\end{equation}
To solve this problem we use the Alternating Direction Multiplier Method (ADMM); see, e.g., \cite{BPNC11} for a recent review. We use ADMM since it allows us to decouple the $\ell^2$ and $\ell^1$ norms as well as the matrices $A$ and $L_\omega$. The latter point is extremely relevant since, as we discuss below, both matrices have exploitable structures, however, they are difficult to exploit together. 

We first reformulate \eqref{eq:model} in an equivalent form
$$
\arg\min_{\x,\y,\w,\z}\left\{\frac{1}{2}\norm{A\x-\bd}^2+\mu\norm{\z}_1+\iota_0(\w),\mbox{ s.t. }\x=\y, \x=\w, \z=L_\omega\y,\right\}
$$
or equivalently
\begin{equation}\label{eq:min_split}
\begin{split}
\arg\min_{\x,\y,\w,\z}&\frac{1}{2}\norm{A\x-\bd}^2+\mu\norm{\z}_1+\iota_0(\w)\\
\mbox{ s.t. }&\begin{bmatrix}
I&O\\
O&I\\
I&O
\end{bmatrix}\begin{bmatrix}
\x\\
\z
\end{bmatrix}-\begin{bmatrix}
I&O\\
L_\omega&O\\
O&I
\end{bmatrix}\begin{bmatrix}
\y\\\w
\end{bmatrix}=\v{0},
\end{split}
\end{equation}
where $O$ and $\v{0}$ denote the zero matrix and the zero vector, respectively, and $\iota_0$ is the indicator function of the nonnegative cone, i.e.,
$$
\iota_0(\x)=\left\{\begin{array}{ll}0&\mbox{if }\x\geq0,\\+\infty&\mbox{otherwise.}\end{array}\right.
$$ 
We can construct the augmented Lagrangian of \eqref{eq:min_split}
\begin{align*}
\mathcal{L}_\rho(\x,\y,\w,\z;\lb)=&\frac{1}{2}\norm{A\x-\bd}^2+\mu\norm{\z}_1+\iota_0(\w)\\&-\sprod{\lb}{\begin{bmatrix}
	I&O\\
	O&I\\
	I&O
	\end{bmatrix}\begin{bmatrix}
	\x\\
	\z
	\end{bmatrix}-\begin{bmatrix}
	I&O\\
	L_\omega&O\\
	O&I
	\end{bmatrix}\begin{bmatrix}
	\y\\\w
	\end{bmatrix}}\\&+\frac{\rho}{2}\norm{\begin{bmatrix}
	I&O\\
	O&I\\
	I&O
	\end{bmatrix}\begin{bmatrix}
	\x\\
	\z
	\end{bmatrix}-\begin{bmatrix}
	I&O\\
	L_\omega&O\\
	O&I
	\end{bmatrix}\begin{bmatrix}
	\y\\\w
	\end{bmatrix}}^2,
\end{align*}
where $\rho>0$ is a fixed parameter and $\lb\in\R^{3N}$ is the Lagrangian multiplier. Applying ADMM we get the iterations
$$
\left\{\begin{array}{l}
\vspace*{2mm}\begin{bmatrix}
\x^{(k+1)}\\\z^{(k+1)}
\end{bmatrix}=\arg\min_{\x,\z}\mathcal{L}_\rho(\x,\y^{(k)},\w^{(k)},\z;\lb^{(k)}),\\
\begin{bmatrix}
\y^{(k+1)}\\\w^{(k+1)}
\end{bmatrix}=\arg\min_{\y,\w}\mathcal{L}_\rho(\x^{(k+1)},\y,\w,\z^{(k+1)};\lb^{(k)}),
\\
\lb^{(k+1)}=\lb^{(k)}+\rho\left(\begin{bmatrix}
I&O\\
O&I\\
I&O
\end{bmatrix}\begin{bmatrix}
\x^{(k+1)}\\
\z^{(k+1)}
\end{bmatrix}-\begin{bmatrix}
I&O\\
L_\omega&O\\
O&I
\end{bmatrix}\begin{bmatrix}
\y^{(k+1)}\\\w^{(k+1)}
\end{bmatrix}\right).
\end{array}
\right.
$$
We can write $\lb=\begin{bmatrix}\lb_1\\\lb_2\\\lb_3\end{bmatrix}$ with $\lb_j\in\R^N$ for $j=1,2,3$. Therefore the minimization problems above decouples and we obtain
$$
\resizebox{\textwidth}{!}{$
\left\{\begin{array}{l}
\x^{(k+1)}=\arg\min_\x \frac{1}{2}\norm{A\x-\bd}^2-\sprod{\begin{bmatrix}\lb_1^{(k)}\\\lb^{(k)}_3\end{bmatrix}}{\begin{bmatrix}\x-\y^{(k)}\\\x-\w^{(k)}\end{bmatrix}}+\frac{\rho}{2}\norm{\begin{bmatrix}\x-\y^{(k)}\\\x-\w^{(k)}	\end{bmatrix}}^2,\\
\z^{(k+1)}=\arg\min_\z \mu\norm{\z}_1-\sprod{\lb_2^{(k)}}{\z-L_\omega\y^{(k)}}+\rho_2\norm{\z-L_\omega\y^{(k)}},\\
\y^{(k+1)}=\arg\min_\y \frac{\rho}{2}\norm{\x^{(k+1)}-\y}^2-\sprod{\begin{bmatrix}\lb_1^{(k)}\\\lb_2^{(k)}\end{bmatrix}}{\begin{bmatrix}\x^{(k+1)}-\y\\L_\omega\z^{(k+1)}-\y\end{bmatrix}}+\frac{\rho}{2}\norm{\begin{bmatrix}\x^{(k+1)}-\y\\L_\omega\z^{(k+1)}-\y\end{bmatrix}}^2,\\
\w^{(k+1)}=\arg\min_\w \frac{\rho}{2}\norm{\x^{(k+1)}-\w}^2-\sprod{\lb_3^{(k)}}{\x^{(k+1)}-\w}+\iota_0(\w).
\end{array}
\right.$}
$$
All of these minimization problems have a closed form for their solution, namely
$$
\left\{\begin{array}{l}
\x^{(k+1)}=(A^TA+2\rho I)^{-1}(A^T\bd+\rho\y^{(k)}-\lb^{(k)}_1+\rho\w^{(k)}-\lb^{(k)}_3),\\
\z^{(k+1)}=S_{\mu/\rho}\left(L_\omega\y^{(k)}-\lb^{(k)}_2/\rho\right),\\
\y^{(k+1)}=(L_\omega^TL_\omega+I)^{-1}(L_\omega^T(\z^{(k+1)}+\lb^{(k)}_2/\rho)+\x^{(k+1)}+\lb^{(k)}_1/\rho),\\
\w^{(k+1)}=\left(\x^{(k+1)}+\lb_3/\rho \right)_+,\\
\lb_1^{(k+1)}=\lb_1^{(k)}+\rho(\x^{(k+1)}-\y^{(k+1)}),\\
\lb_2^{(k+1)}=\lb_2^{(k)}+\rho(\z^{(k+1)}-L_\omega\y^{(k+1)}),\\
\lb_3^{(k+1)}=\lb_3^{(k)}+\rho(\x^{(k+1)}-\w^{(k+1)}),
\end{array}
\right.
$$
where $S_\mu$ denotes the soft-thresholding operator with parameter $\mu$, i.e., $S_\mu(\x)=\ts{sign}(\x)(|\x|-\mu)_+$, where the operations are meant element-wise and $(x)_+=\max\{x,0\}$ is the metric projection into the nonnegative cone. Note that each iteration requires the solution of two linear systems. The linear system involving $A$ can be easily solved using the \texttt{fft} algorithm if periodic BCs are employed; see above. If other BCs are employed, the structure of the matrix, in general, does not allow to use fast transform for the solution of the system. Nevertheless, this system can be solved by an iterative method using a circulant preconditioner. On the other hand, the solution of linear system with the $L_\omega$ matrix can be easily computed using the \texttt{lsqr} algorithm applied to the equivalent least-square problem since the matrix $L_\omega$ is extremely sparse.

We would like to briefly discuss the use of the \texttt{lsqr} method for the solution of the system
\begin{equation}\label{eq:normal}
(L_\omega^TL_\omega+I)\y^{(k+1)}=(L_\omega^T(\z^{(k+1)}+\lb^{(k)}_2/\rho)+\x^{(k+1)}+\lb^{(k)}_1/\rho.
\end{equation}
The linear system of equations \eqref{eq:normal} is equivalent to the least squares problem
\begin{equation}\label{eq:lsqr}
	\begin{split}
\y^{(k+1)}&=\arg\min_{\y}\norm{\begin{bmatrix}
		L_\omega\\I
\end{bmatrix}\y-\begin{bmatrix}
\z^{(k+1)}+\lb^{(k)}_2/\rho\\\x^{(k+1)}+\lb^{(k)}_1/\rho
\end{bmatrix}}^2\\&=\arg\min_{\y}\norm{\widehat{L}\y-\widehat{\v{v}}^{(k)}}^2.
\end{split}
\end{equation}
The \texttt{lsqr} algorithm is an iterative method that determines an approximate solution of \eqref{eq:lsqr} into a Krylov subspace. In particular,
at its $j-$th iteration the \texttt{lsqr} algorithm determines a solution of \eqref{eq:lsqr} in the Krylov subspace $\mathcal{K}_j(\widehat{L}^T\widehat{L},\widehat{L}^T\widehat{\v{v}}^{(k)})$, where
$$
\mathcal{K}_j(\widehat{L}^T\widehat{L},\widehat{L}^T\widehat{\v{v}}^{(k)})=\ts{span}\left\{\widehat{L}^T\widehat{\v{v}}^{(k)},(\widehat{L}^T\widehat{L})\widehat{L}^T\widehat{\v{v}}^{(k)},\ldots,(\widehat{L}^T\widehat{L})^{j-1}\widehat{L}^T\widehat{\v{v}}^{(k)}\right\}.
$$
The \texttt{lsqr} method requires one matrix-vector product with $\widehat{L}$ and one with $\widehat{L}^T$ at each iteration. Therefore, since $\widehat{L}$ is extremely sparse the flop count per iteration is $O(N)$. Moreover \texttt{lsqr} is mathematically equivalent to the \texttt{cg} method applied to \eqref{eq:normal}. However, its implementation is more stable. Nevertheless, the number of iteration required to converge is proportional to $\kappa(\widehat{L})$ which is extremely small; see, e.g., \cite{B96} for a discussion on \texttt{lsqr} and \texttt{cg}.

We summarize our approach in Algorithm~\ref{algo:deblur}.

\begin{algorithm}\caption{Graph Laplacian image deblurring}\label{algo:deblur}
	\begin{algo}
		\INPUT $A\in\R^{N\times N}$, $\bd\in\R^N$, $R>0$, $\sigma>0$, $\rho>0$, $\tau>0$, $K>0$, $\mu>0$
		\OUTPUT $\x\in\R^{n\times n}$
		\STATE Run Algorithm~\ref{algo:1} to compute $L_\omega$
		\STATE $\y^{(0)}=\w^{(0)}=\lb_1^{(0)}=\lb_2^{(0)}=\lb_3^{(0)}=\v{0}$
		\FOR $k=0,\ldots,K$
		\STATE $\x^{(k+1)}=(A^TA+2\rho I)^{-1}(A^T\bd+\rho\y^{(k)}-\lb^{(k)}_1+\rho\w^{(k)}-\lb^{(k)}_3)$
		\STATE $\z^{(k+1)}=S_{\mu/\rho}\left(L_\omega\y^{(k)}-\lb^{(k)}_2/\rho\right)$
		\STATE $\y^{(k+1)}=(L_\omega^TL_\omega+I)^{-1}(L_\omega^T(\z^{(k+1)}+\lb^{(k)}_2/\rho)+\x^{(k+1)}+\lb^{(k)}_1/\rho)$
		\STATE $\w^{(k+1)}=P_0\left(\x^{(k+1)}+\lb_3/\rho \right)$
		\STATE $\lb_1^{(k+1)}=\lb_1^{(k)}+\rho(\x^{(k+1)}-\y^{(k+1)})$
		\STATE $\lb_2^{(k+1)}=\lb_2^{(k)}+\rho(\z^{(k+1)}-L_\omega\y^{(k+1)})$
		\STATE $\lb_3^{(k+1)}=\lb_3^{(k)}+\rho(\x^{(k+1)}-\w^{(k+1)})$
		\IF $k>1\; \& \;\norm{\x^{k+1}-\x^{k}}\leq\tau\norm{\x^{k}}$
		\STATE Exit
		\ENDIF
		\ENDFOR
		\STATE $\x=\x^{(k+1)}$
	\end{algo}
\end{algorithm}

Note that since the functional \eqref{eq:constr_l2l1} is convex we can apply the following classical result on ADMM
\begin{theorem}[see, e.g., Section~3.2 of Boyd et al. \cite{BPNC11}]
With the notation of Algorithm~\ref{algo:deblur} it holds that
\begin{enumerate}[(i)]
	\item $\lim_{k\rightarrow \infty}\norm{\begin{bmatrix}
			I&O\\
			O&I\\
			I&O
		\end{bmatrix}\begin{bmatrix}
			\x^{(k)}\\
			\z^{(k)}
		\end{bmatrix}-\begin{bmatrix}
			I&O\\
			L_\omega&O\\
			O&I
		\end{bmatrix}\begin{bmatrix}
			\y^{(k)}\\\w^{(k)}
	\end{bmatrix}}=0$, i.e., the iterates approach feasibility as $k\rightarrow\infty$;
	\item $\lim_{k\rightarrow \infty}\frac{1}{2}\norm{A\x^{(k)}-\bd}^2+\mu\norm{\z^{(k)}}_1+\iota_0(\w^{(k)})=p^*$, where $p^*$ is the minimum of \eqref{eq:model};
	\item $\lim_{k\rightarrow \infty}\lb_k=\lb^*$, where $\lb^*$ is a dual optimal point, i.e., a saddle point of $L_0$.	
\end{enumerate}
\end{theorem}

\begin{remark}
	ADMM can be slow to converge in certain scenarios. It is not in the scope of this paper to propose a fast algorithm for the solution of \eqref{eq:model}. Rather we wish to show the potentiality of $L_\omega$ as a regularization operator. Nevertheless, it is possible to accelerate the convergence of ADMM by extrapolation methods to improve the convergence rate of ADMM; see, e.g., \cite{BDD20,GOSB14}.
\end{remark}

\section{Numerical Examples}\label{sect:num}
We now report some selected numerical examples to show the performances of the proposed method. We are particularly interested in showing that the graph Laplacian constructed in Algorithm~\ref{algo:deblur} provides better reconstructions than the classical TV approach. 

We compare the results obtained using $L=L_{\rm TV}$ and $L=L_\omega$ in \eqref{eq:constr_l2l1} with the solution $\x^*$ computed in line $9$ of Algorithm~\ref{algo:deblur}. To compute the solution of \eqref{eq:constr_l2l1} with $L=L_{\rm TV}$ we use the algorithm described in \cite{CTY13}. Moreover, we want to show the fill potentiality of the proposed approach. To this aim we construct the operator $L_\omega$ using the exact image $\x_{\rm true}$ and we denote it by $\widetilde{L}_\omega$. Obviously this approach is not feasible in realistic scenarios, nevertheless it allows us to show the full potentiality of the proposed method.  

In all our experiments we set the parameters as specified in Table~\ref{tbl:param}. We compare the considered method in terms of accuracy using the Relative Restoration Error (RRE) computed as
$$
{\rm RRE}(\x)=\frac{\norm{\x-\x_{\rm true}}_2}{\norm{\x_{\rm true}}_2},
$$
the Pick Signal to Noise Ration (PSNR), defined by
$$
{\rm PSNR}(\x)=20\log_{10}\left(\frac{Nm}{\norm{\x_{\rm true}}_2}\right),
$$
where $m$ denotes the maximum value achievable by $\x_{\rm true}$. Moreover, we consider the Structure SIMilarity index (SSIM), constructed in \cite{SSIM}. The definition of the SSIM is extremely involved, here we simply recall that this statistical index measures how structurally similar two images are, in particular, the higher the SSIM the more similar the images are, and its highest achievable value is $1$. 
\begin{table}
	\caption{Setting of the parameters in Algorithm~\ref{algo:deblur}.}
	\label{tbl:param}
	\begin{center}
		\begin{tabular}{lll}
			Paramter&Value&Description\\\hline
			$R$&$10$&Support of the weight function in the Graph\\
			$\sigma$&$10^{-2}$& Variance of the weight function in the Graph\\
			$\rho$&$10^{-1}$& Augmentation parameter in ADMM\\
			$\tau$&$10^{-4}$& Stopping criterion for ADMM\\
			$K$&$3000$&Maximum number of iterations\\
			$\mu$&Hand-tuned&Regularization parameter in \eqref{eq:constr_l2l1}
		\end{tabular}
	\end{center}
\end{table}

\begin{figure}
	\centering
	\begin{minipage}{0.3\textwidth}
		\centering
		\includegraphics[width=\textwidth]{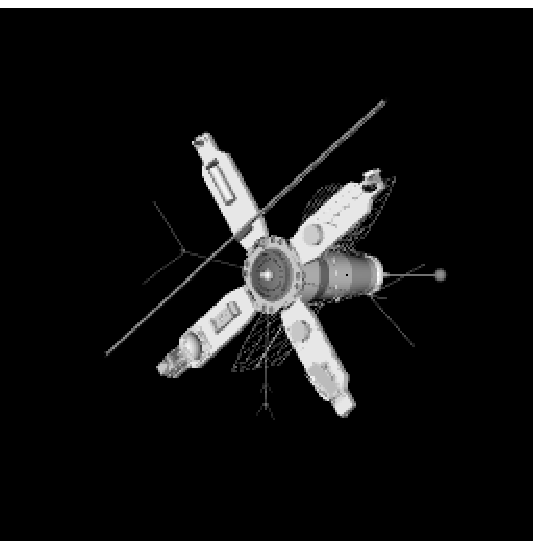}\\(a)
	\end{minipage}
	\begin{minipage}{0.3\textwidth}
		\centering
		\includegraphics[width=\textwidth]{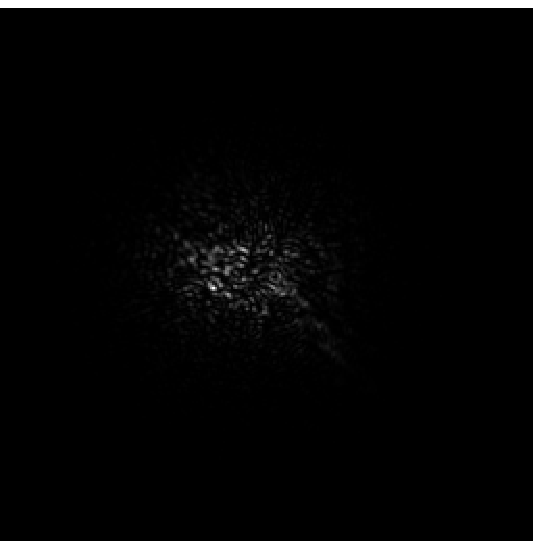}\\(b)
	\end{minipage}
	\begin{minipage}{0.3\textwidth}
		\centering
		\includegraphics[width=\textwidth]{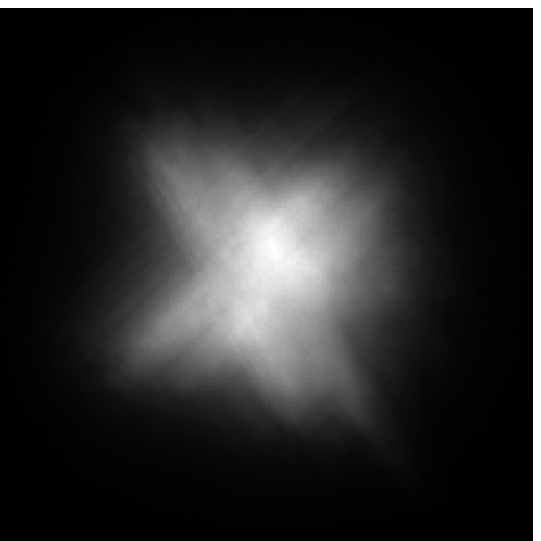}\\(c)
	\end{minipage}
	\caption{Example 1: (a) True image ($256\times 256$ pixels), (b) PSF ($256\times 256$ pixels), (c) Blurred and noisy image ($256\times 256$ pixels with $\delta\approx0.01\norm{\b}_2$).}
	\label{fig:ex1}
\end{figure}

\paragraph{Example 1.} Our first example is the \texttt{atmosphericBlur50} test case of the RestoreTools toolbox \cite{restoreTools}. We report the exact image, the PSF, and the blurred an noisy image in Fig.~\ref{fig:ex1}. The norm of the noise, denoted by $\delta$, that corrupts the data is approximately $1\%$ of the norm of the exact right-hand side $\bd$.

\begin{figure}
	\centering
	\includegraphics[width=0.3\textwidth]{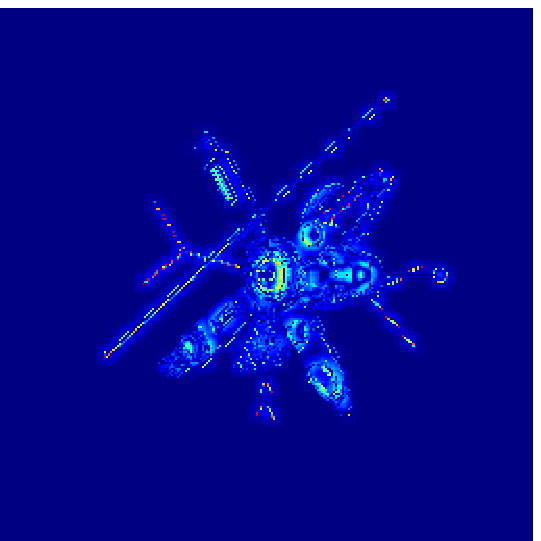}
	\caption{Example 1: Visualization of $\left|\widetilde{L}_\omega\x^\dagger\right|$ in the \texttt{jet} colormap. The color blue represents the $0$s in the image.}
	\label{fig:sparsity}
\end{figure}

\begin{figure}
	\centering
	\begin{minipage}{0.3\textwidth}
		\centering
		\includegraphics[width=\textwidth]{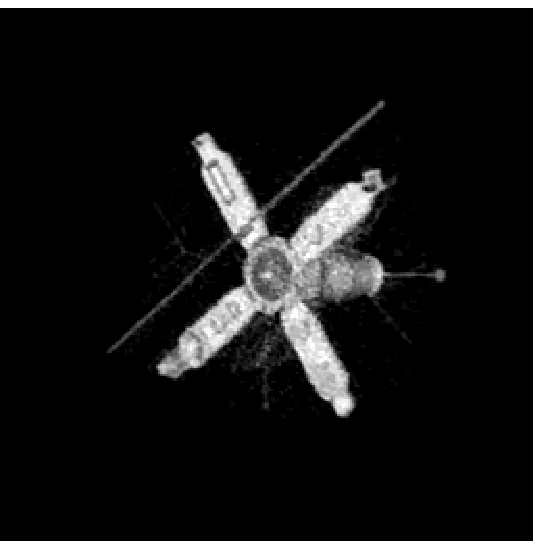}\\(a)
	\end{minipage}
	\begin{minipage}{0.3\textwidth}
		\centering
		\includegraphics[width=\textwidth]{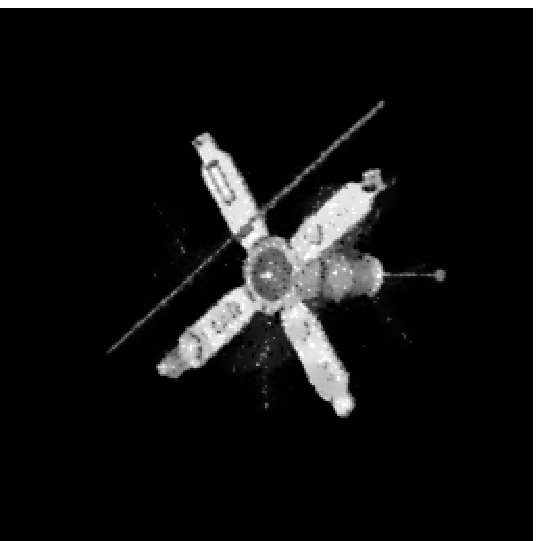}\\(b)
	\end{minipage}
	\begin{minipage}{0.3\textwidth}
		\centering
		\includegraphics[width=\textwidth]{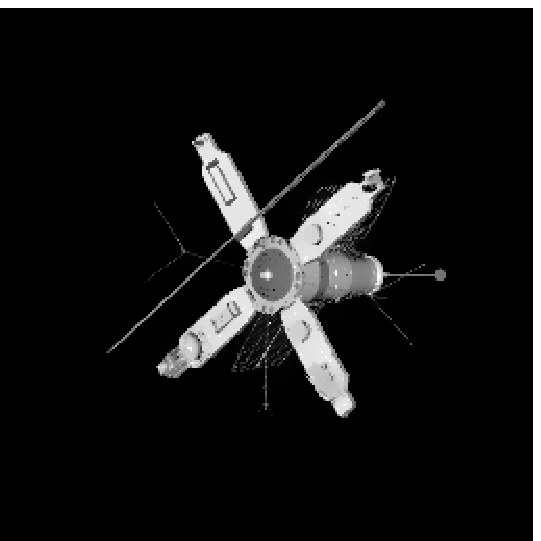}\\(c)
	\end{minipage}
	\caption{Example 1 reconstructions: (a) $\ell^2-\ell^1$ with $L=L_{\rm TV}$, (b) $\ell^2-\ell^1$ with $L=L_{\omega}$, (c) $\ell^2-\ell^1$ with $L=\widetilde{L}_{\omega}$.}
	\label{fig:ex1_rec}
\end{figure}

We report the obtained results with the considered methods in Table~\ref{tbl:res}. We can observe that $\ell^2-\ell^1$ methods provide much more accurate results than the classical Tikhonov method, especially in terms of SSIM. The reconstruction obtained with $L=\widetilde{L}_\omega$, i.e., using the graph related to the exact image, is extremely accurate. To validate our model we show in Fig.~\ref{fig:sparsity} a visualization of $\left|\widetilde{L}_\omega\x^\dagger\right|$. We can observe that it is extremely sparse. Therefore, we expect our model to provide accurate reconstructions. However, this approach is not feasible in real scenarios. Nevertheless, we can also observe that using $L=L_{\omega}$ improves the quality of the restoration with respect to the classic TV. This is confirmed by the visual inspection of the reconstructions in Fig.~\ref{fig:ex1_rec}. Comparing the reconstructions obtained by the $\ell^2-\ell^1$ methods we can observe that the choice $L=L_{\rm TV}$ leads to more noisy reconstructions than the ones obtained with the graph Laplacian.

\begin{figure}
	\centering
	\begin{minipage}{0.3\textwidth}
		\centering
		\includegraphics[width=\textwidth]{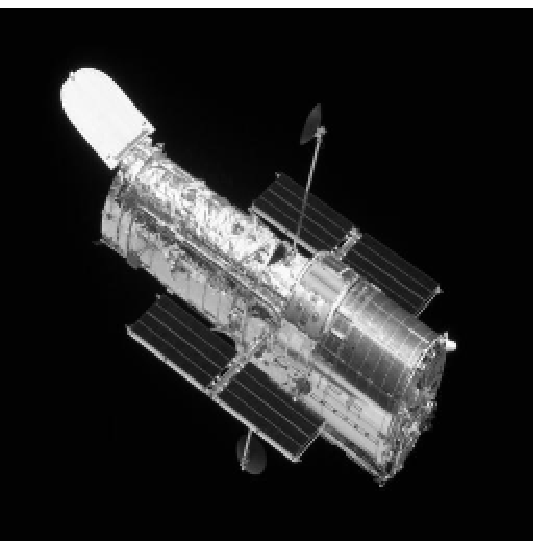}\\(a)
	\end{minipage}
	\begin{minipage}{0.3\textwidth}
		\centering
		\includegraphics[width=\textwidth]{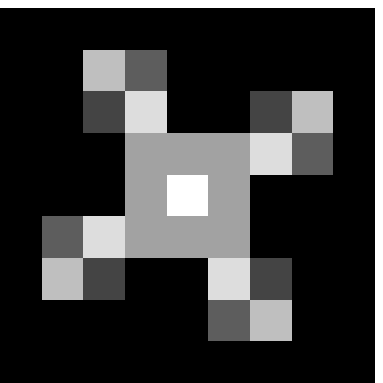}\\(b)
	\end{minipage}
	\begin{minipage}{0.3\textwidth}
		\centering
		\includegraphics[width=\textwidth]{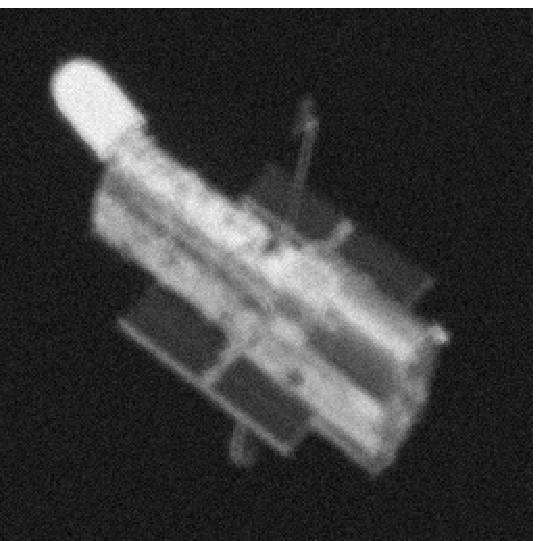}\\(c)
	\end{minipage}
	\caption{Example 2: (a) True image ($256\times 256$ pixels), (b) PSF ($9\times 9$ pixels), (c) Blurred and noisy image ($256\times 256$ pixels with $\delta=0.1\norm{\b}_2$).}
	\label{fig:ex2}
\end{figure}

\paragraph{Example 2.}
For our second example we consider the Hubble image in Fig.~\ref{fig:ex2}(a) and we blur it with the PSF in Fig.~\ref{fig:ex2}(b). We then add white Gaussian noise such that $\delta=0.1\norm{\b}_2$ obtaining the blurred and noisy image in Fig.~\ref{fig:ex2}(c).

\begin{figure}
	\centering
	\begin{minipage}{0.3\textwidth}
		\centering
		\includegraphics[width=\textwidth]{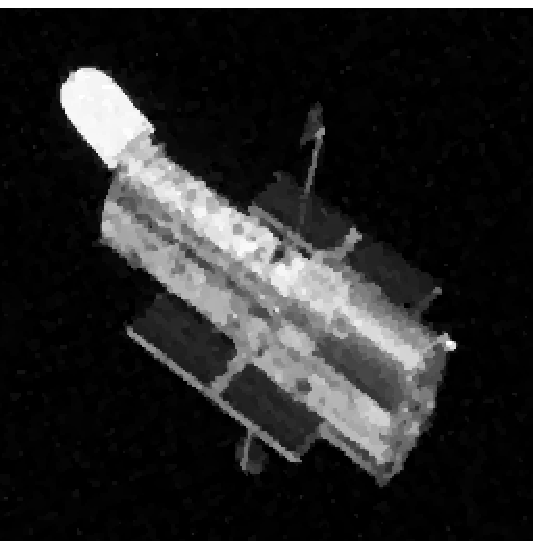}\\(a)
	\end{minipage}
	\begin{minipage}{0.3\textwidth}
		\centering
		\includegraphics[width=\textwidth]{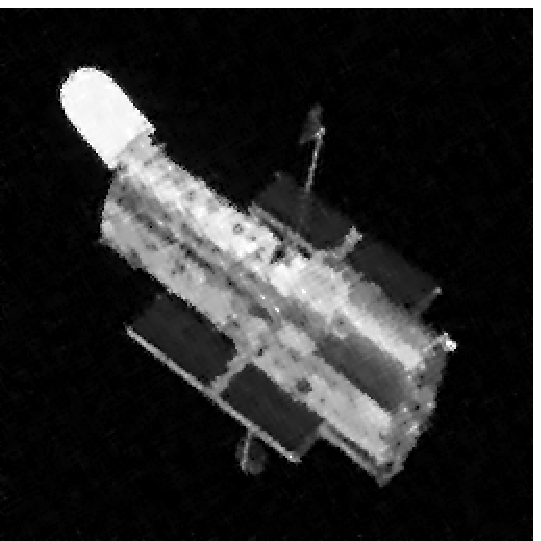}\\(b)
	\end{minipage}
	\begin{minipage}{0.3\textwidth}
		\centering
		\includegraphics[width=\textwidth]{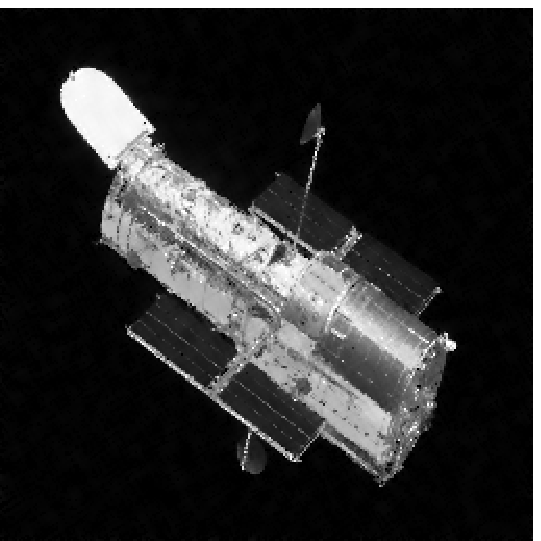}\\(c)
	\end{minipage}
	\caption{Example 2 reconstructions: (a) $\ell^2-\ell^1$ with $L=L_{\rm TV}$, (b) $\ell^2-\ell^1$ with $L=L_{\omega}$, (c) $\ell^2-\ell^1$ with $L=\widetilde{L}_{\omega}$.}
	\label{fig:ex2_rec}
\end{figure}

\begin{figure}
	\centering
	\begin{minipage}{0.3\textwidth}
		\centering
		\includegraphics[width=\textwidth]{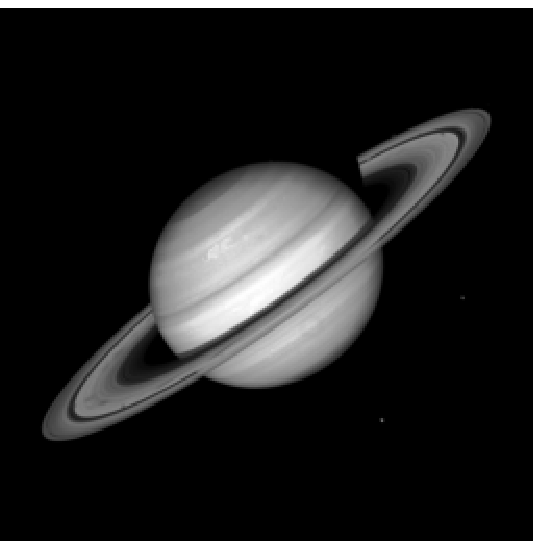}\\(a)
	\end{minipage}
	\begin{minipage}{0.3\textwidth}
		\centering
		\includegraphics[width=\textwidth]{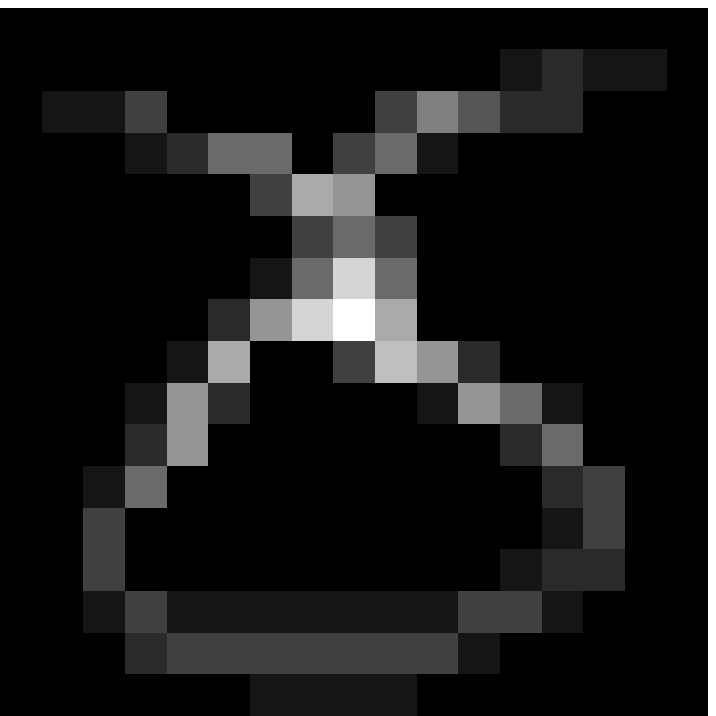}\\(b)
	\end{minipage}
	\begin{minipage}{0.3\textwidth}
		\centering
		\includegraphics[width=\textwidth]{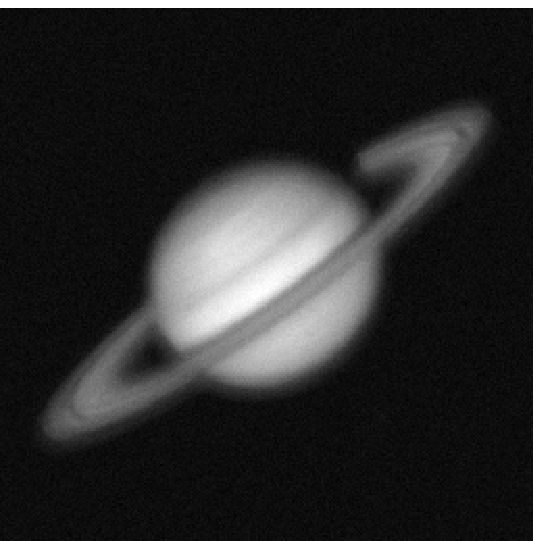}\\(c)
	\end{minipage}
	\caption{Example 3: (a) True image ($256\times 256$ pixels), (b) PSF ($17\times 17$ pixels), (c) Blurred and noisy image ($256\times 256$ pixels with $\delta=0.05\norm{\b}_2$).}
	\label{fig:ex3}
\end{figure}

We compute approximate solutions with the considered algorithms and report the obtained RRE, PSNR, and SSIM in Table~\ref{tbl:res}. We can observe that our proposal provides the best reconstruction both in terms of RRE (and, therefore of PSNR) and SSIM. This is confirmed by the visual inspection of the reconstructions in Fig.~\ref{fig:ex2_rec}. We would like to stress that, similarly as the previous example, the unconstrained Tikhonov method computes extremely noisy reconstructions.
\paragraph{Example 3.}
For our third example we consider the Saturn image in Fig.~\ref{fig:ex3}(a). We blur it with a non-symmetric PSF (see Fig.~\ref{fig:ex3}(b)) and add $5\%$ of white Gaussian noise, i.e., $\delta=0.05\norm{\b}_2$ obtaining the image in Fig.~\ref{fig:ex3}(c).

\begin{figure}
	\centering
	\begin{minipage}{0.3\textwidth}
		\centering
		\includegraphics[width=\textwidth]{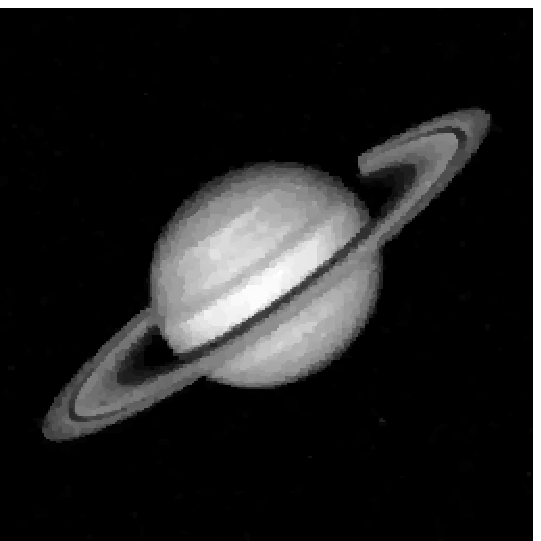}\\(a)
	\end{minipage}
	\begin{minipage}{0.3\textwidth}
		\centering
		\includegraphics[width=\textwidth]{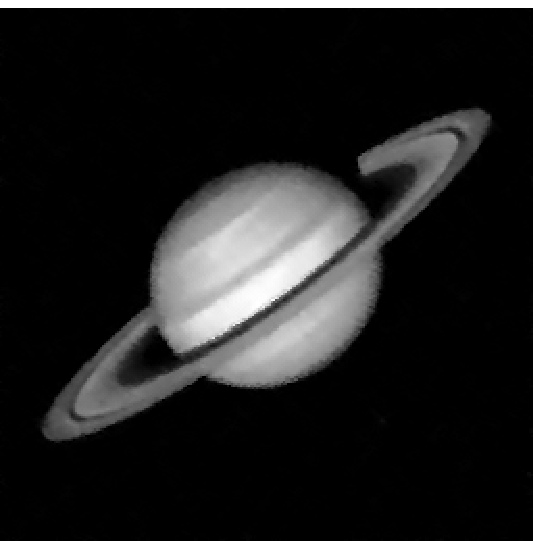}\\(b)
	\end{minipage}
	\begin{minipage}{0.3\textwidth}
		\centering
		\includegraphics[width=\textwidth]{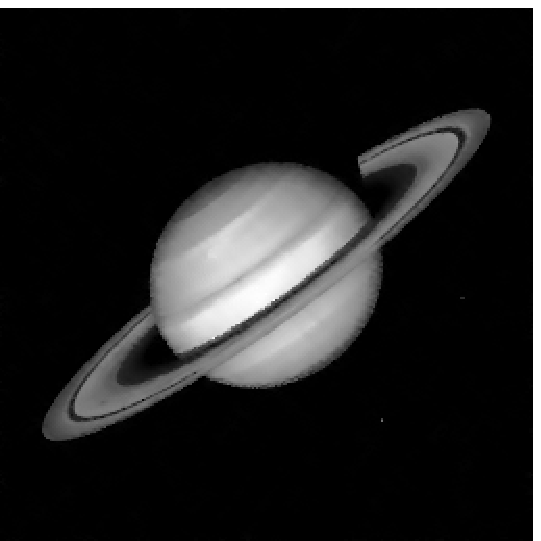}\\(c)
	\end{minipage}
	\caption{Example 3 reconstructions: (a) $\ell^2-\ell^1$ with $L=L_{\rm TV}$, (b) $\ell^2-\ell^1$ with $L=L_{\omega}$, (c) $\ell^2-\ell^1$ with $L=\widetilde{L}_{\omega}$.}
	\label{fig:ex3_rec}
\end{figure}

\begin{figure}
	\centering
	\begin{minipage}{0.3\textwidth}
		\centering
		\includegraphics[width=\textwidth]{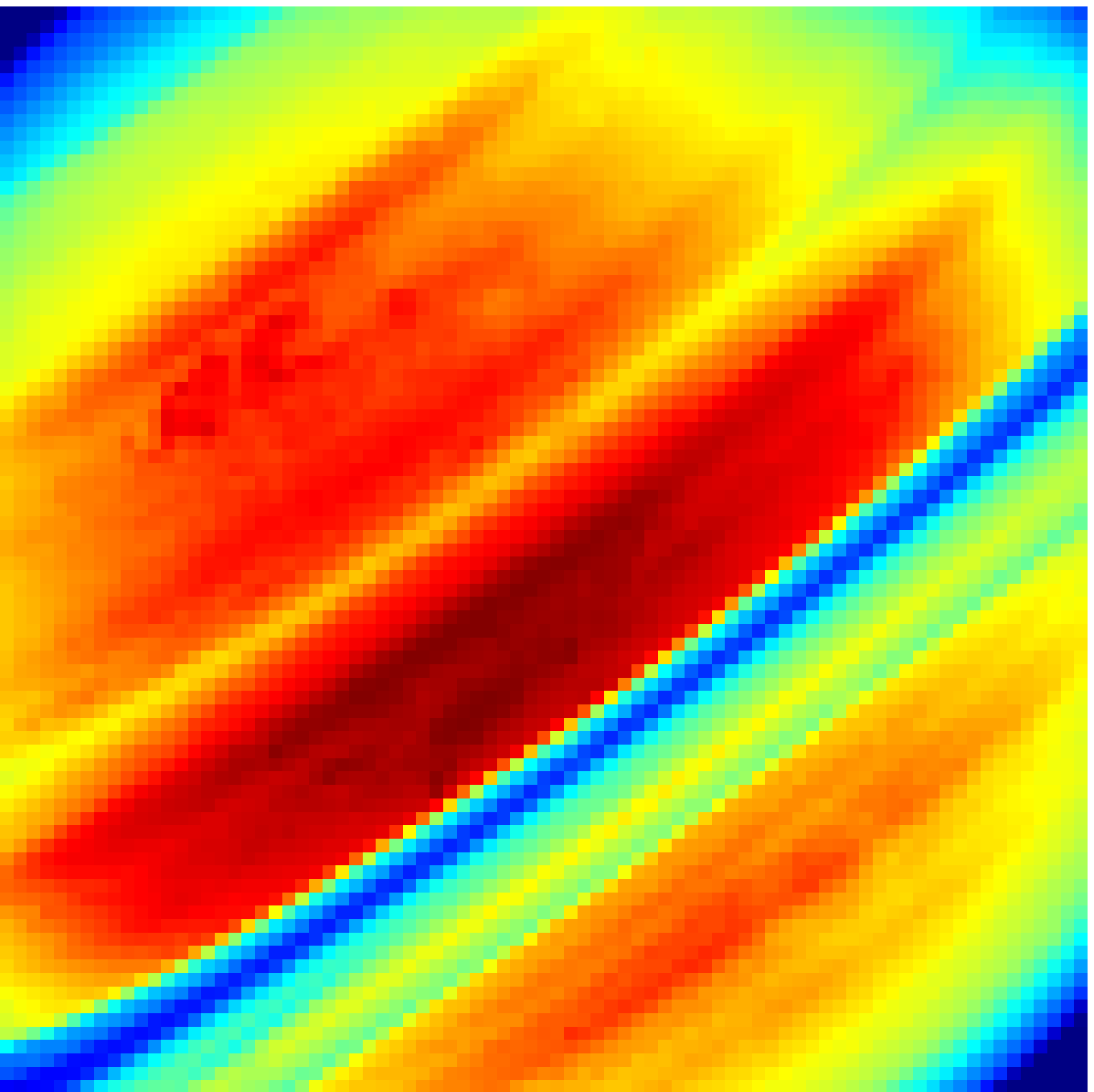}\\(a)
	\end{minipage}
	\begin{minipage}{0.3\textwidth}
		\centering
		\includegraphics[width=\textwidth]{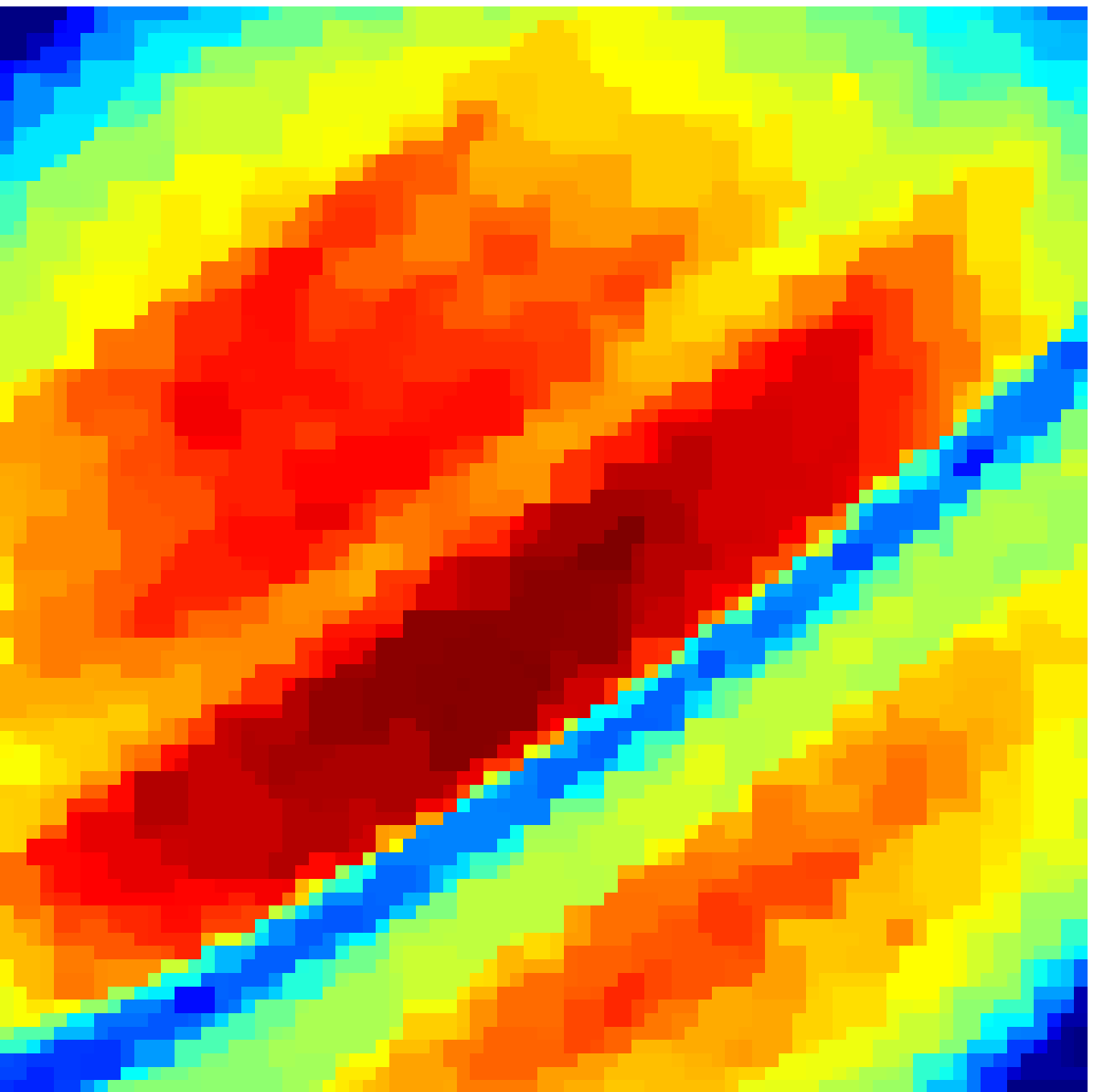}\\(b)
	\end{minipage}
	\begin{minipage}{0.3\textwidth}
		\centering
		\includegraphics[width=\textwidth]{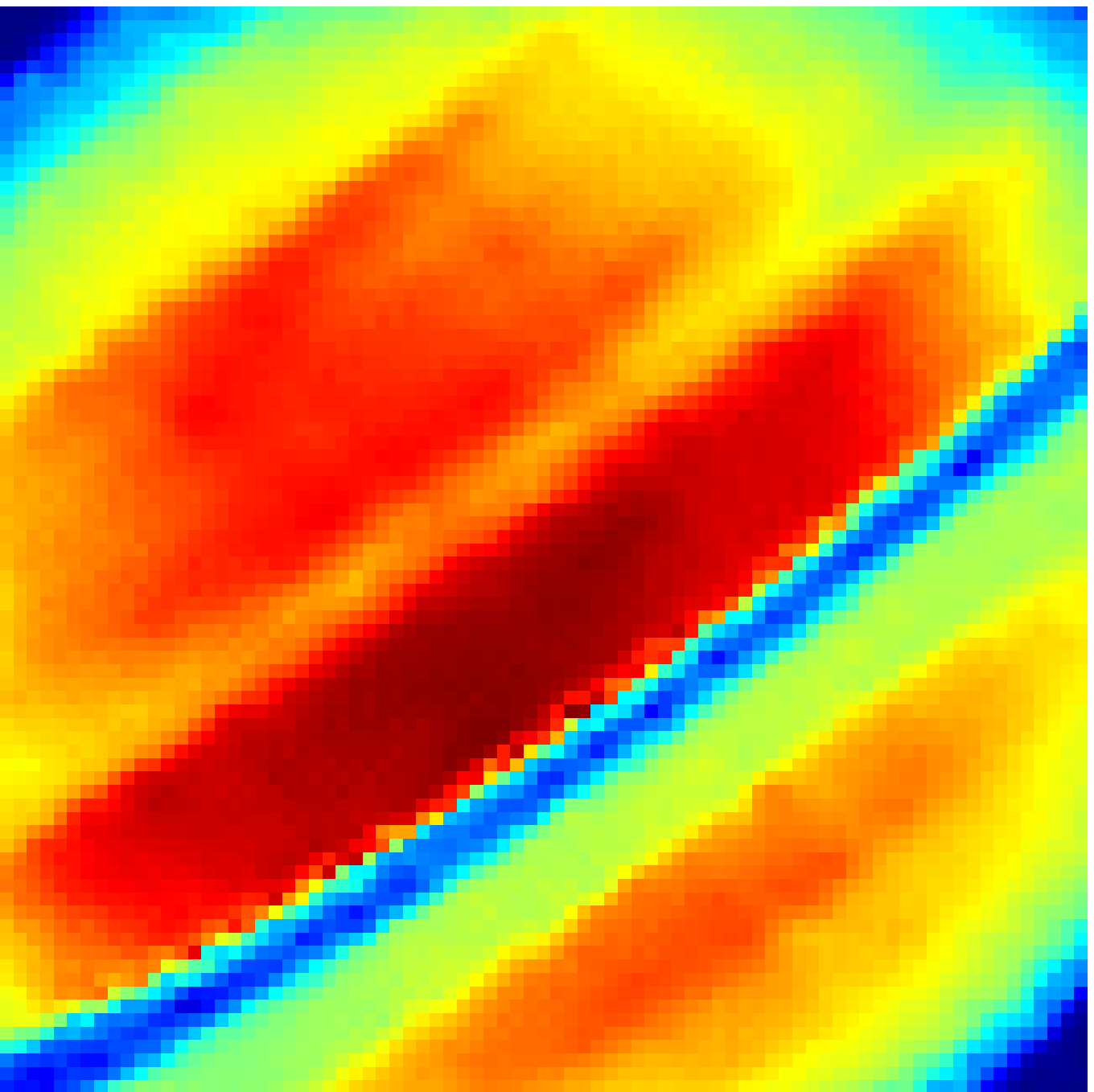}\\(c)
	\end{minipage}
	\caption{Example 3 blow ups of the exact solution and of two reconstructions in the \texttt{jet} colormap: (a) True solution, (b) $\ell^2-\ell^1$ with $L=L_{\rm TV}$, (c) $\ell^2-\ell^1$ with $L=L_{\omega}$.}
	\label{fig:ex3_rec_zoom}
\end{figure}

We report in Table~\ref{tbl:res} the obtained results with the considered algorithms. We can observe that our proposal provides a very accurate reconstruction in terms of RRE and PSNR. However, the SSIM of the computed solution is slightly lower than the one obtained with the standard TV regularization. In Fig.~\ref{fig:ex3_rec} we report all the computed solution. We would like to observe that the reconstruction obtained with the classical TV regularization presents a very heavy stair-case effect while the approximate solution obtained with our proposal avoids this issue. To show this we propose in Fig.~\ref{fig:ex3_rec_zoom} blow-ups of the central part of the image of the exact solution and of the reconstructions obtained by TV regularization and our approach. We can observe that the TV reconstruction presents a very heavy stair-case effect that is avoided by our proposal. The solution computed by our method manages to remain extremely sharp, while avoiding any stair-casing. 

\begin{figure}
	\centering
	\begin{minipage}{0.3\textwidth}
		\centering
		\includegraphics[width=\textwidth]{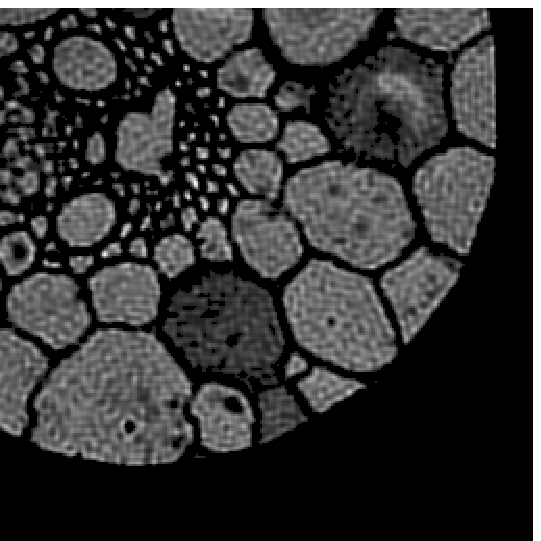}\\(a)
	\end{minipage}
	\begin{minipage}{0.3\textwidth}
		\centering
		\includegraphics[width=\textwidth]{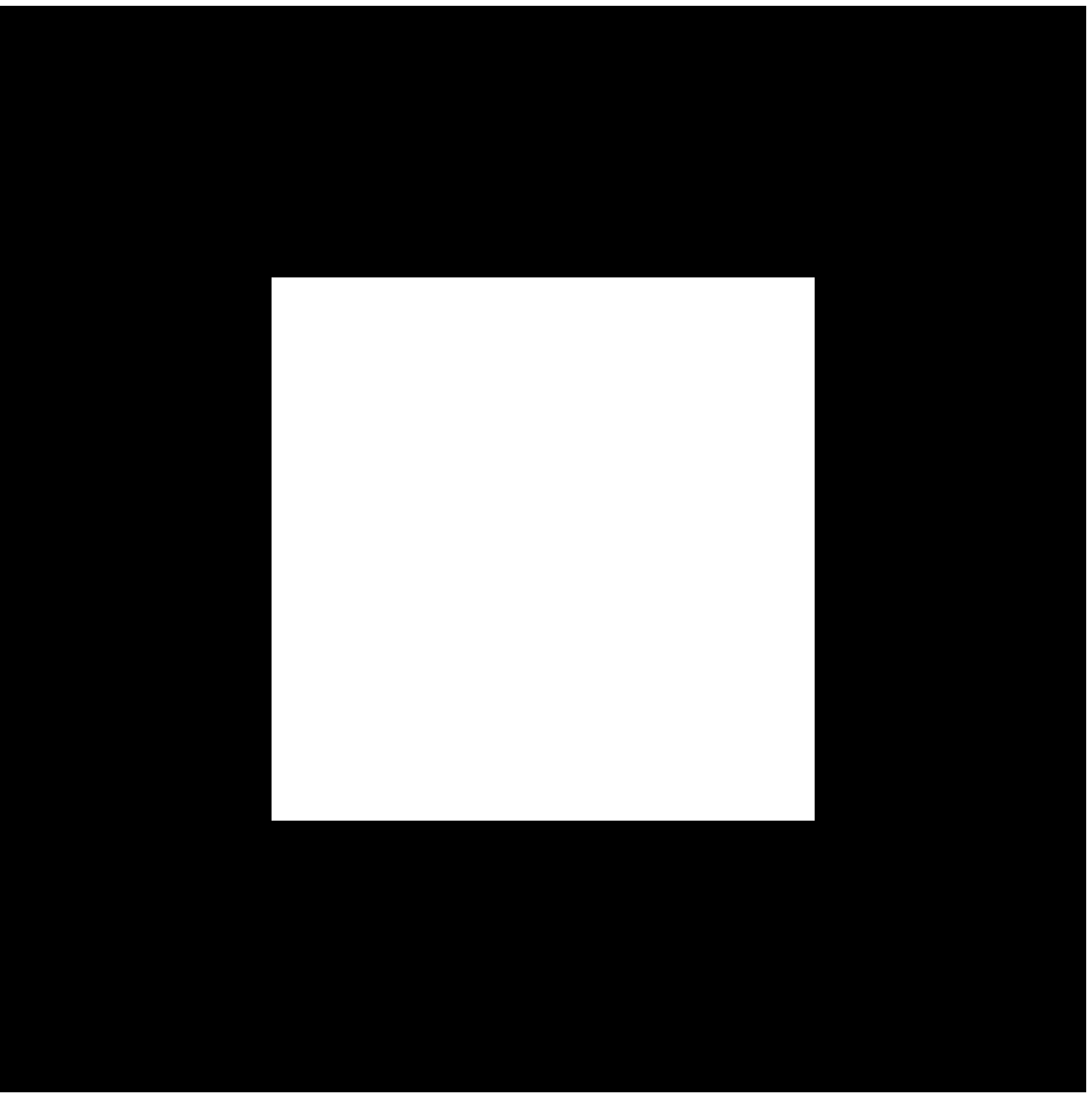}\\(b)
	\end{minipage}
	\begin{minipage}{0.3\textwidth}
		\centering
		\includegraphics[width=\textwidth]{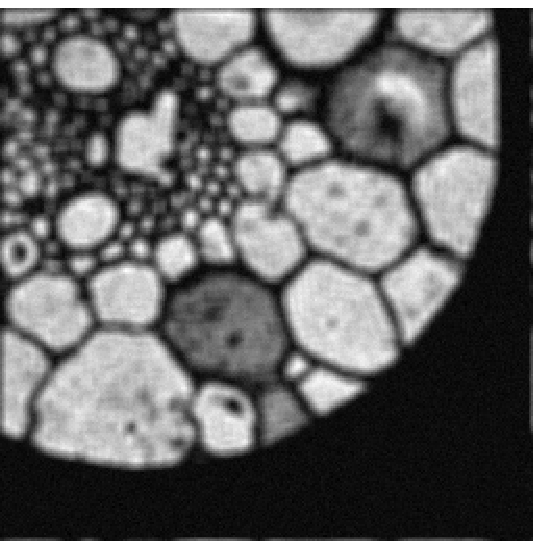}\\(c)
	\end{minipage}
	\caption{Example 4: (a) True image ($256\times 256$ pixels), (b) PSF ($12\times 12$ pixels), (c) Blurred and noisy image ($256\times 256$ pixels with $\delta=0.03\norm{\b}_2$).}
	\label{fig:ex0}
\end{figure}

\begin{figure}
	\centering
	\begin{minipage}{0.3\textwidth}
		\centering
		\includegraphics[width=\textwidth]{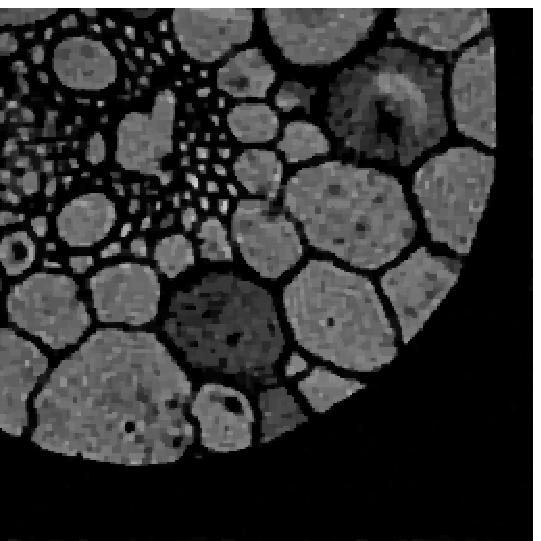}\\(a)
	\end{minipage}
	\begin{minipage}{0.3\textwidth}
		\centering
		\includegraphics[width=\textwidth]{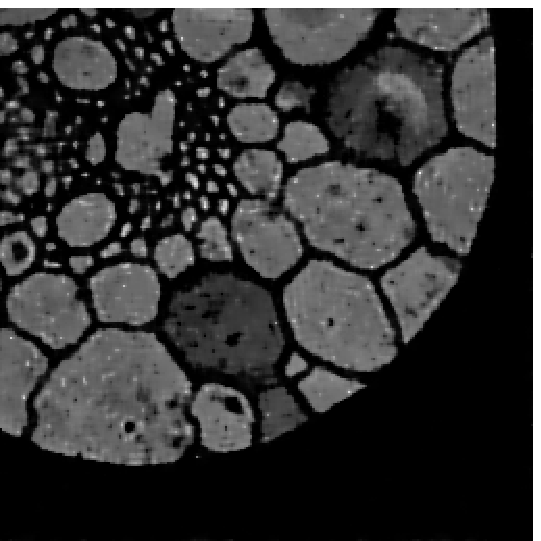}\\(b)
	\end{minipage}
	\begin{minipage}{0.3\textwidth}
		\centering
		\includegraphics[width=\textwidth]{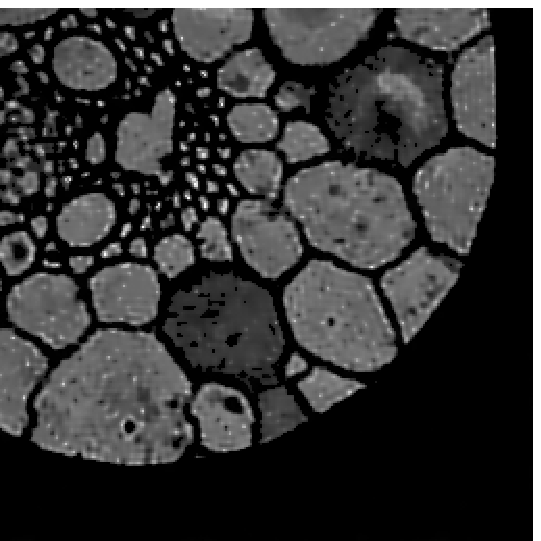}\\(c)
	\end{minipage}
	\caption{Example 4 reconstructions: (a) $\ell^2-\ell^1$ with $L=L_{\rm TV}$, (b) $\ell^2-\ell^1$ with $L=L_{\omega}$, (c) $\ell^2-\ell^1$ with $L=\widetilde{L}_{\omega}$.}
	\label{fig:ex0_rec}
\end{figure}

\paragraph{Example 4.}We consider the exact image in Fig.~\ref{fig:ex0}(a) and we blur it with an average PSF (see Fig.~\ref{fig:ex0}(b)), we add white Gaussian noise so that $\norm{\b-\bd}=0.03\norm{\b}$ obtaining Fig.~\ref{fig:ex0}(c). We construct the graph Laplacian related to the exact solution $\widetilde{L}_\omega$ and deblur the blurred and noisy image with Tikhonov (with $\mu=\mu_{\rm GCV}$) and by minimizing \eqref{eq:constr_l2l1} with both $L=L_{\rm TV}$ and $L=\widetilde{L}_\omega$. We report the results obtained in Table~\ref{tbl:res}. We can observe that the proposed algorithm with the perfect choice of $\widetilde{L}_\omega$ largely outperforms the other approaches furnishing a very accurate reconstruction of the proposed image. Moreover, the choice of $L=L_\omega$, which we recall that it does not require any a priori information on the exact solution, is still more accurate than the classical TV.  This is confirmed by the visual inspection of the reconstructions in Fig.~\ref{fig:ex0_rec}.


\begin{table}
	\begin{center}
		\caption{Comparison of the RRE, PSNR, and SSIM.}
		\label{tbl:res}
		\begin{tabular}{l|l|lll}
			Example&Method&RRE&PSNR&SSIM\\\hline
			\multirow{4}{*}{Example 1}&Tikhonov&$0.22299$&$26.663$&$0.55512$\\
			&$\ell^2-\ell^1$ with $L=L_{\rm TV}$&$0.19152$&$27.984$&$0.92623$\\
			&$\ell^2-\ell^1$ with $L=L_{\omega}$&${0.17763}$&${28.638}$&${0.93971}$\\
			&$\ell^2-\ell^1$ with $L=\widetilde{L}_{\omega}$&$0.083333$&$35.212$&$0.98129$\\\hline
			\multirow{4}{*}{Example 2}&Tikhonov&$0.17352$&$25.735$&$0.55241$\\
			&$\ell^2-\ell^1$ with $L=L_{\rm TV}$&$0.15492$&$26.720$&$0.80458$\\
			&$\ell^2-\ell^1$ with $L=L_{\omega}$&$0.14968$&$27.019$&$0.81256$\\
			&$\ell^2-\ell^1$ with $L=\widetilde{L}_{\omega}$&$0.10096$&$30.439$&$0.89943$\\\hline
			\multirow{4}{*}{Example 3}&Tikhonov&$0.080283$&$33.715$&$0.72254$\\
			&$\ell^2-\ell^1$ with $L=L_{\rm TV}$&$0.068917$&$35.041$&$0.94873$\\
			&$\ell^2-\ell^1$ with $L=L_{\omega}$&$0.060094$&$36.231$&$0.94809$\\
			&$\ell^2-\ell^1$ with $L=\widetilde{L}_{\omega}$&$0.040382$&$39.684$&$0.96489$\\\hline
			\multirow{4}{*}{Example 4}&Tikhonov&$0.16236$&$27.160$&$0.73224$\\
			&$\ell^2-\ell^1$ with $L=L_{\rm TV}$&$0.15299$&$27.686$&$0.86899$\\
			&$\ell^2-\ell^1$ with $L=L_{\omega}$&$0.14716$&$28.024$&$0.85765$\\
			&$\ell^2-\ell^1$ with $L=\widetilde{L}_{\omega}$&$0.085936$&$32.695$&$0.93887$
		\end{tabular}
	\end{center}
\end{table}
\section{Conclusions}\label{sect:concl}
In this paper we have proposed a new regularization operator for $\ell^2-\ell^1$ minimization. The construction of this operator is automatic and extremely cheap to perform. We have shown that the proposed method outperforms the classical TV approach. Matter of future research include the application of the proposed method to more general inverse problems as well as the integration of the considered method with the $\ell^p-\ell^q$ minimization proposed in \cite{HLMSR17,LMSR15,BR19,BRXX,BRXXb,BPR20} or with iterative regularization methods like, e.g., Linearized Bregman splitting \cite{COS09,COS09b,BPR18,CDBH16,BPRXXb} and Iterated Tikhonov with general penalty terms \cite{BBDS15,BDR16,B17,Bianchi2017}. Another line of future research is the construction of more sophisticated graphs $\omega$ which can better exploit the structure of the given image itself. Such constructions may stem from a PDEs approach; see, e.g., \cite{Gilboa2009,Lou2010,Adriani2020}.

\section*{Acknowledgments}
A.B., and M.D. are members of the GNCS-INdAM group. D.B. is member of the GNAMPA-INdAM group. A.B. research is partially supported by the Regione Autonoma della Sardegna research project ``Algorithms and Models for Imaging Science [AMIS]'' (RASSR57257, intervento finanziato con risorse FSC 2014-2020 - Patto per lo Sviluppo della Regione Sardegna). 

\bibliography{biblio}
\bibliographystyle{spmpsci}
\end{document}